\documentclass{amsart}

\usepackage{amsmath}
\usepackage{amsthm}
\usepackage{amsfonts}
\usepackage{amssymb}

\renewcommand\P{{\mathbf{P}}}

\newcommand\R{{\mathbf{R}}}
\newcommand\Z{{\mathbf{Z}}}
\newcommand\C{{\mathbf{C}}}
\newcommand\E{{\mathbf{E}}}

\newcommand\I{{\mathbf{I}}}

\newcommand\N{{\mathcal{N}}}

\newcommand\tr{{\operatorname{tr}}}
\renewcommand\th{{\operatorname{th}}}

\newcommand\Incomp{{\operatorname{Incomp}}}
\newcommand\Comp{{\operatorname{Comp}}}

\newcommand\dist{{\operatorname{dist}}}

\renewcommand\Re{{\operatorname{Re}}}
\renewcommand\Im{{\operatorname{Im}}}

\renewcommand\a{{x}}

\newcommand\eps{\varepsilon}

\theoremstyle{plain}
  \newtheorem{theorem}[subsection]{Theorem}

  \newtheorem{lemma}[subsection]{Lemma}
  \newtheorem{corollary}[subsection]{Corollary}
   
\theoremstyle{remark}
  \newtheorem{remark}[subsection]{Remark}
  
\theoremstyle{definition}
  \newtheorem{definition}[subsection]{Definition}
\include{psfig}
\begin{document}

\title[ ]
{Smooth analysis of the condition number and the least singular value}

\author{Terence Tao}
\address{Department of Mathematics, UCLA, Los Angeles CA 90095-1555}
\email{tao@math.ucla.edu}
\thanks{T. Tao is supported by a grant from the MacArthur Foundation.}
\author{  Van  Vu}
\address{Department of Mathematics, Rutgers, Piscataway, NJ 08854}
\email{vanvu@math.rutgers.edu}
\thanks{V.   Vu  is supported by  NSF  Grant DMS-0901216 and AFOSAR-FA-9550-09-1-0167.}


\begin{abstract}
Let $\a$ be a complex  random variable with mean zero and bounded
variance. Let $N_{n}$ be the random matrix of size $n$ whose entries
are iid copies of $\a$ and $M$ be a fixed matrix of the same size. The goal of this paper is to give  a
general  estimate for the condition number and  least singular value of the matrix
$M + N_{n}$, generalizing an earlier result of Spielman and Teng for the case when
$\a$ is gaussian.

Our investigation reveals an  interesting fact that  the ``core'' matrix $M$
does play a role on  tail bounds for the least singular value
of $M+N_{n} $.  This does  not occur in Spielman-Teng studies when $\a$ is gaussian.
 Consequently, our general estimate involves the norm $\|M\|$.
 In the special case when $\|M\|$ is relatively small,  this estimate is nearly  optimal and  extends
or refines existing  results.
\end{abstract}
\maketitle

\section {Introduction}

Let $M$ be an $n \times n$ matrix and $s_{1 } (M) \ge \dots \ge s_{n } (M)$ its singular values.
The condition number of $A$, as defined by numerical analysts,  is

$$\kappa (M) := s_{1} (M)/ s_{n} (M) = \|M\| \|M^{-1} \|. $$

This parameter is of fundamental importance in numerical linear algebra and
related areas, such as linear programming. In particular, the value

$$L(M):= \log \kappa (M) $$ measures  the (worst case) lost of precision the equation
$Mx=b$ can exhibit \cite{Smale, BT}.

The problem of  understanding  the typical  behavior of $\kappa(M)$ and $L(M)$
when the matrix $M$ is random has  a long history. This  was first
raised by von Neuman and Goldstine in their study of numerical inversion of large
matrices \cite{vonG}. Several years later, the problem was restated in a survey of Smale \cite{Smale} on the efficiency of algorithm of anaylsis. One of Smale's motivations was to understand the efficiency of the simplex algorithm
in linear programming.  The problem is also at the core of  Demmel's plan about the investigation  of
the probability that a numerical analysis problem is difficult \cite{Demmel} (see also \cite{Ren} for a work that inspires this investigation).

To make the problem precise, the  most critical issue  is to choose
a probability distribution for $M$.
A convenient  model has been  random matrices with independent
gaussian entries (either real of complex).  An essential feature of this model is that here
the joint distribution of the eigenvalues can be written down precisely

\begin{equation} \label{jointreal} (Real \,\,\, Gaussian)\,\,\,
c_1(n) \prod_{1 \le i < j \le n} |\lambda_i -\lambda_j| \exp(- \sum_{i=1}^n \lambda_i^2/2 ).
\end{equation}

\begin{equation} \label{jointcomplex} (Complex \,\,\, Gaussian)\,\,\,
c_2 (n) \prod_{1 \le i < j \le n} |\lambda_i -\lambda_j|^2 \exp(- \sum_{i=1}^n \lambda_i^2/2 ).
\end{equation}

Here $c_1(n), c_2(n)$ are normalization factors whose explicit formulae can be seen in, for example, \cite{Mehta}.

Most questions about  the spectrum   of these random matrices can then be answered by
estimating a properly defined  integral with respect to these measures.
Many advanced techniques have been worked
out to serve this purpose (see, for instance \cite{Mehta}).
 In particular, the condition number is well understood, thanks to works of
 Kostlan, Oceanu \cite{Smale, Kos}, Edelman \cite{Edel} and many others
 (see Section \ref{section:previous}).

The gaussian model, however,  has serious shortcomings.
As pointed out by many researchers (see, for example \cite{BCL0, ST}), the gaussian model
does not reflex the arbitrariness of the input. Let us consider, for example,
 a random matrix with independent real gaussian entries.
By sharp concentration results, one can show
that  the fraction of entries with absolute values at most 1, is,
with overwhelming probability, close to the absolute constant
$\frac{1}{\sqrt {2 \pi}} \int_{-1}^1 \exp(-t^2/2) dt$.
Many classes of matrices that
occur in practice just simply do not posses this property. This problem persists even
when one replaces gaussian by another fixed distribution, such as Bernoulli.

About 10 years ago, Spielman and Teng \cite{ST, ST1}, motivated by Demmel's plan and
the problem of understanding the efficiency of the simplex algorithm
 proposed a new, exciting  distribution.
  Spielman and Teng observed that while
the ideal input maybe a fixed matrix $M$, it is likely that the computer will work with a
perturbation $M+N$, where $N$ is a random matrix representing random noise. Thus, it raised  the issue of studying the distribution of the  condition number of $M+N$. This problem is at  the heart of the so-called
Spielman-Teng {\it smooth } analysis. (See \cite{ST, ST1} for a more detailed discussion and
\cite{BCL0, BCL1, BCL2, SST, DST} for many related
works on this topics.)  Notice that the special case $M=0$ corresponds to
the setting considered in the previous paragraphs.

Spielman-Teng model nicely addresses the problem about the arbitrariness of the inputs, as in this model every matrix generates a probability space of its own.
In their papers, Spielman and Teng   considered mostly
 gaussian noise  (in some cases they also considered other continuous distributions such as uniform on $[-1,1]$). However, in  the digital world, randomness  often {\it does not}
 has gaussian nature. To start with, all of real
  data are finite.
In fact, in many problems (particularly those in integer programming)  all entries of the matrix   are integers. The random errors made by the degital
devices (for example, sometime a bit gets flipped) are
obviously of discrete nature.
 In other problems, for example those in engineering, the data may
contain measurements where
it would be  natural to assume gaussian errors. On the other hand,
data are usually strongly truncated. For example, if an entry of our matrix represents
 the mass of an object, then we expect  to see  a number like $12.679$ (say, tons), rather than
$12.6792347043641259$. Thus, instead of the gaussian distribution, we (and/or our computers) often work with a discrete distribution, whose support is relatively small and does not depend
on the size of the matrix. (A good toy example is random Bernoulli matrix, whose entries takes values $\pm 1$ with probability half.)  This leads us to the following question

\vskip2mm

\vskip2mm
{\bf Question.} (Smooth analysis of the condition number)  Estimate  the condition number
of a random matrix $M_{n} := M+N_{n}$, where $M$ is a fixed matrix of size $n$, and $N_{n} $  a general random matrix ?

\vskip2mm

The goal of  this paper is  to
investigate this question, where, as a generalization of Spielman-Teng model, we
 think of $N_{n}$ as a matrix with independent random entries which (instead as being gaussian)  have arbitrary distributions.  Our main result will show that with high probability, $M_{n} $ is well-conditioned.
 This result could be useful in further studies of smooth analysis in linear programming.
 The Spielman-Teng smooth analysis of the simplex algorithm \cite{ST, ST1} was done with gaussian noise. It is a natural and (from the practical point of view) important question to
 repeat this analysis with discrete noise (such as Bernoulli). This question was posed by Spielman to
 the authors few years ago. The paper \cite{ST} also contains a specific conjecture
 on the least singular value of random Bernoulli matrix.

 In connection, we should mention here  a recent series of papers  by  Burgisser, Cucker and Lotz \cite{BCL0, BCL1, BCL2}, which discussed the smooth analysis of condition number  under a  somewhat different setting (they considered the notion of {\it conic } condition number and a different kind of randomness).

Before stating mathematical results, let us describe our notations.
We use the usual asymptotic notation $X = O(Y)$ to
denote the estimate $|X| \leq CY$ for some constant $C
>0$ (independent of $n$); $X=\Omega(Y)$ to denote the estimate $X \geq cY$ for some $c > 0$
independent of $n$, and $X = \Theta(Y)$ to denote the estimates
$X=O(Y)$ and $X=\Omega(Y)$ holding simultaneously. In some cases, we
write $X \ll Y$ instead of $X=O(Y)$ and $X \gg Y$ instead of
$X=\Omega (Y)$. Notations such as $X= O_{\a,b} (Y)$ or $X \ll_{a,b}
(Y)$ mean that the hidden constant in $O $ or $\ll$ depend on
previously defined constants $a$ and $b$.
We use $o(1)$ to denote any quantity that goes to zero as $n \to
\infty$. $X=o(Y)$ means that $X/Y=o(1)$.

\vskip2mm

Recall that

$$\kappa (M) := s_{1} (M)/ s_{n} (M) = \|M\| \|M^{-1} \|. $$

Since $\|M\|^{2} \ge \sum_{ij} |m_{ij} |^{2} / n$ (where $m_{ij}$ denote the entries of $M$) it is expected that
$\|M\| = n^{\Omega (1)}$. Following the literature, we say that
$M$ is well-conditioned (or well-posed) if $\kappa (M) = n^{O(1)} $ or (equivalently)
$L(M)= O(\log n)$.

By the triangle inequality,

$$\|M\| -\|N_{n}\| \le \| M+N_{n}\| \le \|M\| + \|N_{n}\| . $$

Under very general assumptions, the random matrix $N_{n}$ satisfies
$\| N_{n}\| = n^{O(1)} $ with overwhelming probability (see many estimates in Section \ref{section:newresults}). Thus, in order to guarantee that
$\|M+N_{n}\|$ is well-conditioned (with high probability), it is natural to assume that

\begin{equation}  \label{eqn:assumption0} \|M\| = n^{O(1)} . \end{equation}

This is not only a natural, but fairly safe assumption to make (with respect to the applicability of
our studies). Most large
 matrices in practice satisfy this assumption, as their  entries are usually not too large
 compared to their sizes.

 Our main result shows that under this assumption
and a very general assumption on the entries of $N_{n}$, the matrix
$M+N_{n}$ is well-conditioned, with high probability. This result extends
and bridges several existing results in the literature (see next two sections).

Notice that under assumption \eqref{eqn:assumption0}, if we want to show that
$M+N_{n}$ is typically well-conditioned, it suffices to show that

$$\| (M+N_{n})^{-1} \| = s_{n} (M+N_{n})^{-1} = n^{O(1) } $$ with high probability. Thus, we will formulate
most results  in a form of a tail bound for the least singular value of $M+N_{n}$. The typical form will be

$$\P (s_{n} (M+N_{n}) \le n^{-B}) \le n^{-A} $$ where $A,B $ are positive constants
and $A $ increases  with $B$. The relation between $A$ and $B$ is of importance and will be
discussed in length.

\section{Previous results} \label{section:previous}

Let us first discuss the gaussian case.  Improving results of Kostlan and Oceanu \cite{Smale},
Edelman \cite{Edel} computed the limiting
 distribution of $\sqrt n s_{n} (N_{n}) $ when  $N_{n} $ is gaussian.
His result implies

\begin{theorem} \label{theorem:Edel} There is a constant $C >0$ such that the following holds.
 Let $\a$ be the real gaussian random variable with mean zero and
variance one, let $N_{n}$ be the random matrix whose entries are iid
copies of $\a$. Then  for any constant $t>0$
$$\P (s_{n}(N_{n}) \le t) \le  n^{1/2} t. $$
\end{theorem}

Concerning the more general model $M+N_{n}$,  Sankar, Spielman and Teng proved
\cite{SST}

\begin{theorem} \label{theorem:SST} There is a constant $C>0$ such
that the following holds. Let $\a$ be the real gaussian random
variable with mean zero and variance one, let $N_{n}$ be the random
matrix whose entries are iid copies of $\a$, and let $M$ be an
arbitrary fixed matrix. Let $M_{n}:= M+N_{n}$. Then for any $t
>0$

$$\P(s _n (M_n) \le t) \le C n^{1/2} t. $$
\end{theorem}

Once we give up the gaussian assumption,  the study of the least
singular value $s_n$ becomes much harder  (in particular for
discrete distributions such as Bernoulli, in which $\a = \pm 1$ with
equal probability $1/2$).  For example, it is  already
non-trivial to prove that the least singular value of a  random
Bernoulli matrix is positive with probability $1-o(1)$. This was
first done by Koml\'os in 1967 \cite{Kom}, but good quantitative
lower bounds were not available until recently. In a series of
papers, Tao-Vu and Rudelson-Vershynin addressed this question
\cite{TVsing, TVcir, Rud, RV} and proved a lower bound of the form
$n^{-\Theta(1)}$ for $s_{n}$ with high probability.

We say that $\a$ is {\it subgaussian} if there is a constant $B >0$
such that $$\P (|\a| \ge t) \le 2 \exp(-t^2/ B^2)$$ for all $t >0$.
The smallest $B$ is called the {\it subgaussian moment} of $\a$.
The following is a corollary of a more general theorem by Rudelson and Vershynin
\cite[Theorem 1.2]{RV}

\begin{theorem} \label{theorem:RV} Let $\a$ be a subgaussian random variable with zero mean, variance
one and subgaussian moment $B$ and $A$ be an arbitrary positive
constant.
 Let $N_n$ be the random matrix whose entries
are iid copies of $\a$. Then there is a positive constant $C$
(depending on $B$) such that  for any $t \ge n^{-A}$ we have
$$\P( s_n (N_n) \le t) \le C n^{1/2}t . $$
\end{theorem}

We again turn to the general model $M+N_{n}$. In \cite{TVcir}, the present authors proved

\begin{theorem}\label{cirmain}\cite[Theorem 2.1]{TVcir} Let $\a$ be a random variable
 with non-zero variance. Then for any constants $A , C> 0$ there exists a constant
 $B > 0$ (depending on $A, C$, $\a$) such that the following holds.
  Let  $N_n$ be the random matrix whose entries are iid copies
of $\a$, and let $M$ be any deterministic $n \times n$ matrix with
norm $\|M\| \le n^C$.  Then

$$\P( s_n (M+N_n) \le n^{-B}) \le n^{-A}.$$

\end{theorem}

Notice that this theorem requires very little about  the variable
$\a$. It does not need to be sub-gaussian nor even has bounded moments. All we ask is that the variance is bounded from zero, which basically means $\a$ is indeed ``random''.
Thus, it guarantees the well-conditionness of $M+N_{n} $ in a very general setting.

The weakness of this theorem is that the dependence of $B$ on $A$
and $C$, while explicit, is too
generous. The main result of this paper, Theorem \ref{theorem:main},
will improve this dependence significantly and provide a common
extension of Theorem \ref{cirmain} and Theorem \ref{theorem:RV}.

\section{Main result}  \label{section:newresults}

As already pointed out, an important point is the
relation between the constants $A,B$ in a bound of the form

$$\P(s_{n} (M+N_{n} ) \le n^{- B} ) \le n^{-A} . $$

In  Theorem \ref{theorem:SST}, we have a  simple (and optimal) relation
$B= A+1/2$.  It is  natural to conjecture that this relation holds for other, non-gaussian, models of random matrices. In fact, this conjecture  was our starting point of this study.
Quite surprisingly, it turns out not to be the case.

\begin{theorem} \label{theorem:construction} There are positive constants $c_1$ and $c_2$ such that the following holds.
Let $N_n$ be the $n \times n$ random Bernoulli matrix with $n$ even.
For any $L\ge n$, there is an $n \times n$ deterministic matrix $M$
such that $\|M\|= L$ and
$$ \P( s_n(M+N_n) \le c_1 \frac{n}{L} ) \ge c_2 n^{-1/2}.$$
\end{theorem}

The assumption $n$ is even is for convenience and can easily be removed by replacing the Bernoulli matrix by a random matrix whose entries take values $0, \pm 1$ with probability $1/3$ (say).
Notice that if  $L= n^{D}$ for some constant $D$ then we have the lower bound

$$ \P( s_n(M+N_n) \le c_{1 }n^{-D+1 }) \ge c_2 n^{-1/2},$$

which shows that one cannot expect Theorem \ref{theorem:SST} to hold in general and
that the norm of $M$ should play a role in tail bounds of the least singular value.

The main result of this paper is the following.

\begin{theorem} \label{theorem:main} Let $\a$ be a random variable with mean zero and bounded
second moment,
 and let $\gamma \geq 1/2$, $A \ge 0$ be constants. Then
 there is a constant $c$ depending on $\a, \gamma, A$ such
 that the following holds.
 Let $N_{n}$ be the random matrix of size $n$ whose entries are iid copies of $\a$,
  $M$ be a deterministic
   matrix satisfying $\|M\| \le  n^\gamma$, and let $M_{n}:=
   M+N_{n}$. Then

$$\P (s_n (M_n) \leq n^{-(2A+1)\gamma}) \le c \Big( n^{-A+o(1)} + \P( \|N_n\| \geq n^\gamma ) \Big).$$

\end{theorem}

Note that this theorem only assumes bounded second moment on $\a$. The assumption that
the entries of $N_{n}$ are iid is for convenience. A slightly weaker result would hold if one omit this assumption.

\begin{corollary} \label{cor:main} Let $\a$ be a random variable with mean zero and bounded
second moment,
 and let $\gamma \geq 1/2$, $A \ge 0$ be constants. Then
 there is a constant $c_2$ depending on $\a, \gamma, A$ such
 that the following holds.
 Let $N_{n}$ be the random matrix of size $n$ whose entries are iid copies of $\a$,
  $M$ be a deterministic
   matrix satisfying $\|M\| \le  n^\gamma$, and let $M_{n}:=
   M+N_{n}$. Then

$$\P (\kappa (M_n) \ge 2n^{(2A+2)\gamma}) \le c \Big( n^{-A+o(1)} + \P( \|N_n\| \geq n^\gamma ) \Big).$$

\end{corollary}

\begin{proof}
Since $\kappa (M_{n} )= s_{1} (M_{n} )/ s_{n } (M_{n}) $, it follows that  if
$\kappa (M_n) \ge n^{(2A+2)\gamma}$, then at least one of the two events
 $s_{n} (M_{n} ) \le n^{-(2A+1)\gamma}$ and $s_{1} (M_{n} ) \ge 2n^{\gamma} $  holds.
 On the other hand, $$s_{1} (M_{n} ) \le s_{1 } (M) + s_{1 } (N_{n})= \|M\| + \|N_{n } \| \le n^{\gamma} +\|N_{n}\|. $$ The claim follows.
\end{proof}

In the rest of this section, we deduce a few corollaries
and connect them with the existing results.

First, consider the special case when $\a$ is subgaussian. In this case, it is well-known that one can
have a strong bound on $\P( \|N_n\| \geq n^\gamma )$ thanks to the
following theorem (see \cite{RV} for references)

\begin{theorem} \label{theorem:largest1} Let $B$ be a positive
constant. There are positive constants $C_1, C_2$ depending on $B$
such that the following holds. Let $\a$ be a subgaussian random
variable with zero mean, variance one and subgaussian moment $B$
 and  $N_n$ be the random matrix whose entries
are iid copies of $\a$.  Then

$$ \P( \| N_{n} \| \geq C_1 n^{1/2} ) \le \exp(-C_2 n).$$

If one replaces the subgaussian condition by the weaker condition
that $\a$ has forth moment bounded  $B$, then one has a weaker
conclusion that

$$\E (\|N_{n}\| ) \le C_1 n^{1/2} . $$
\end{theorem}

From Theorem \ref{theorem:main}  and Theorem \ref{theorem:largest1}
we see that

\begin{corollary} \label{cor:main1}  Let $A$ and $\gamma$ be arbitrary positive constants. Let
$\a$ be a subgaussian random variable with zero mean and variance
one
 and  $N_n$ be the random matrix whose entries
are iid copies of $\a$. Let $M$ be a deterministic matrix such that
$\| M\| \le n^{\gamma}$ and set $M_n=M+ N_n$. Then

\begin{equation}\label{mm2}
 \P( s_n (M_n) \le (n^{1/2} + \|M\|)^{-2A-1} ) \le  n^{-A+o(1)}.
\end{equation}
\end{corollary}

In the case $\| M \| = O(n^{1/2})$ (which of course includes the $M=0$ special case), \eqref{mm2} implies

\begin{corollary} \label{cor:main1-RV}  Let $A$  be arbitrary positive constant. Let
$\a$ be a subgaussian random variable with zero mean and variance
one
 and  $N_n$ be the random matrix whose entries
are iid copies of $\a$. Let $M$ be a deterministic matrix such that
$\| M\| = O(n^{1/2})$ and set $M_n=M+ N_n$. Then

\begin{equation} \label{eqn:RV1} \P( s_n (M_n) \le n^{-A-1/2}  ) \le  n^{-A+o(1)}.
\end{equation}
\end{corollary}

Up to a loss of magnitude $n^{o(1)}$, this matches Theorem
\ref{theorem:RV}, which treated the base case $M=0$.

If we assume bounded fourth moment instead of subgaussian, we can use
the second half of Theorem \ref{theorem:largest1} to deduce

\begin{corollary} \label{cor:main2}  Let $\a$ be a random variable with zero mean, variance one and
bounded forth moment moment
 and  $N_n$ be the random matrix whose entries
are iid copies of $\a$. Let $M$ be a deterministic matrix such that
$\| M\| = n^{O(1)}$ and set $M_n =M+ N_n$. Then

\begin{equation}\label{mm3}
 \P( s_n (M_n) \le  (n^{1/2} + \|M\|)^{-1+o(1)} ) =o(1).
 \end{equation}
\end{corollary}

In the case $\| M\| = O(n^{1/2})$, this implies that almost surely
$s_n (M_n) \ge n^{-1/2 +o(1)}$. For the special case $M=0$, this
matches (again up to the $o(1)$ term) Theorem \cite[Theorem
1.1]{RV}.

 Let us now take a look at the influence of $\| M \|$ on
 the bound. Obviously, there is  a gap between \eqref{mm2} and Theorem \ref{theorem:construction}.
On the other hand, by setting $A=1/2$, $L= n^{\gamma}$ and assuming
that $\P (\|N_n \| \ge n^{\gamma})$ is negligible (i.e.,
super-polynomially small in $n$), we can deduce from Theorem
\ref{theorem:main} that

$$\P (s_n (M_n) \leq c_1 L^{-2} ) \le  c_2 n^{-1/2+o(1)}.$$

This, together with Theorem \ref{theorem:construction}, suggests that
the influence of $\| M\|$ in  $s_n (M_n)$ is of polynomial type.

In the next discussion, let us
normalize and   assume that $\a$ has  variance one. One can deduce a
bound on $\| N_n \|$ from the simple computation

$$ \E \| N_n \|^2 \leq \E  \,\, \tr N_n N_n^* =  n^2 .$$

By Chebyshev's inequality we thus have
$$ \P( \|N_n\| \geq n^{1+A/2} ) \le n^{-A}$$
for all $A \geq 0$.

Applying Theorem \ref{theorem:main} we obtain

\begin{corollary} \label{cor:main3} Let $\a$ be a random variable
with mean zero and variance one and $N_n$ be the random matrix whose
entries are iid copies of $\a$. Then for any constant $A \ge 0$

$$ \P( s_n(N_n) \le n^{-1-\frac{5}{2}A - A^2} ) \le n^{-A+o(1)}. $$

In  particular, $s_n(N_n) \ge n^{-1-o(1)}$ almost surely.
\end{corollary}

It is clear that one can obtain better bounds for $s_n$, provided
better estimates on $\| N_n \|$. The idea of using Chebyshev's
inequality is very crude (we just like to give an example) and there
are more sophisticated tools. One can, for instance, use higher moments.
The expectation of a $k$-th moment can be expressed a sum of many terms, each correspond to a certain closed walk of length $k$ on the complete graph of $n$ vertices (see \cite{FK, Vnorm}).
 If the higher moments of
$N_{n}$ (while not bounded) do not increase too fast with $n$, then the main contribution in
the expectation of the $k$th  moment still come from terms which correspond to walks using each edge of the graph either 0 and 2 times. The expectation of such a term involves only the second moment of
the entries in $N_{n}$. The reader may want to work this out as an exercise.

One can  also use the following  nice estimate of Seginer \cite{Seg}

$$\E \|N_{n} \| =O( \E \max_{1\le i \le n} \sqrt {\sum_{j=1}^{n} \a_{ij}^{2} }  +  \E \max_{1\le j \le n} \sqrt {\sum_{i=1}^{n} \a_{ij}^{2} }). $$

The rest of the paper is organized as follows. In the next section, we prove
Theorem \ref{theorem:construction}. The remaining sections are devoted for the proof of
Theorem \ref{theorem:main}.
This proof combines several tools that have been developed in recent years.
It starts with an  $\epsilon$-net argument (in the spirit of  those used  in \cite{TVsing, Rud,
TVcir, RV}. Two important technical ingredients are Theorem \ref{theorem:TV} from \cite{TVcir} and
Lemma \ref{lemma:RV} from \cite{RV}.

\section{Theorem \ref{theorem:construction}: The influence of $M$}

 Let $M'$ be the $n-1 \times n$ matrix obtained by concatenating
the matrix $LI_{n-1}$ with an all $L$ column, where $L$ is a large
number (we will set $L \ge n$). The $n \times n$ matrix $M$ is obtained from $M'$ by adding
to it a (first) all zero row; thus
$$ M = \begin{pmatrix}
0 & 0 & \ldots & 0 & 0 \\
L & 0 & \ldots & 0 & L \\
0 & L & \ldots & 0 & L \\
\vdots & \vdots & \ddots & \vdots & \vdots \\
0 & 0 & \ldots & L & L
\end{pmatrix}.$$
It is easy to see that

$$\|M\| =\Theta (L). $$

Now consider $M_n:= M+ N_n$ where the entries of $N_n$ are iid
Bernoulli random variables.
$$\P( s_n (M_n) \ll n^{1/4}L^{-1/2}) \gg n^{-1/2}. $$

Let $M'_n$ be the (random) $(n-1) \times n$ matrix formed by the
last $n-1$ rows of $M_n$. Let $v \in \R^n$ be a unit normal vector of the $n-1$ rows of $M'_n$. By replacing $v$ with $-v$ if necessary we may write $v$ in the
form
$$v= \left(\frac{1}{\sqrt n} +a_1,
\frac{1}{\sqrt n}+ a_2, \dots , \frac{1}{\sqrt n} +a_{n-1},
\frac{-1}{\sqrt n} + a_n\right),$$
where $\frac{-1}{\sqrt n} + a_n \le 0$.

Let $\xi_i$ be iid Bernoulli random variables. Multiplying $v$ with the first row of $M'_n$, we have
\begin{align*}
0&= (L+\xi_1)(\frac{1}{\sqrt n} +a_1) +  (L+\xi_n)(-\frac{1}{\sqrt n}+a_n) \\
&= L(a_1+a_n) + \frac{1}{\sqrt n} \Big((\xi_1-\xi_n) + \xi_1a_1 + \xi_na_n \Big).
\end{align*}

Since $|a_i| = O(1)$, it follows that $|a_1+a_n| =O(\frac{1}{L})$.
Repeating the argument with all other rows, we conclude that
$|a_i+a_n| = O(\frac{1}{L})$ for all $1\le i \le n-1$.

Since $v$ has unit norm, we also have
$$1 = \|v\|^2= \sum_{i=1}^{n-1} \left(\frac{1}{\sqrt n} + a_i\right)^2 +
\left(\frac{-1}{\sqrt n} + a_n\right)^2, $$
which implies that
$$ \frac{2}{\sqrt n} (a_1+ \dots + a_{n-1} -a_n) + \sum_{i=1}^n
a_i^2 =0. $$

This, together with the fact that $|a_i+a_n| = O(\frac{1}{L})$ and
all $1\le i \le n-1$, yields

$$na_n^2- 2n a_n (\frac{1}{\sqrt n} + \frac{1}{L})  = O(\frac{ \sqrt n}{L} +\frac{1}{L^2} ) . $$

Since $-\frac{1}{\sqrt n}+a_n \le 0$ and $L \ge n$, it is easy to
show from here that $|a_n|= O(\frac{1}{L})$. It follows that $|a_i|
= O(\frac{1}{L})$ for all $1\le i \le n$.

Now consider

$$\|M_n v \|= \left|\sum_{i=1}^{n-1} (\frac{1}{\sqrt n}+a_i) \xi_i +
(-\frac{1}{\sqrt n}+a_n) \xi_n\right|. $$

Since $n$ is even, with probability $\Theta (\frac{1}{\sqrt n})$,
$\xi_1+ \dots +\xi_{n-1} -\xi_n=0$, and in this case

$$\|M_n v \| = \left| \sum_{i=1}^n a_i \xi_i \right|= O\left(\frac{n}{L}\right), $$
as desired.

\section{Controlled moment} \label{section:generaltheorems}

It is convenient to establish some more quantitative control on $\a$.  We recall the following notion
from \cite{TVcir}.

\begin{definition}[Controlled second moment]  Let $\kappa \geq 1$.  A complex random variable $\a$ is said to have \emph{$\kappa$-controlled second moment} if one has the upper bound
$$ \E |\a|^2 \leq \kappa$$
(in particular, $|\E \a| \leq \kappa^{1/2}$), and the lower bound
\begin{equation}\label{eyot}
 \E \Re( z \a - w )^2 \I(|\a| \leq \kappa) \geq \frac{1}{\kappa} \Re(z)^2
\end{equation}
for all complex numbers $z,w$.
\end{definition}

{\it Example} The Bernoulli random variable ($\P(\a=+1) = \P(\a=-1) = 1/2$) has $1$-controlled second moment.  The condition \eqref{eyot} asserts in particular that $\a$ has variance at least $\frac{1}{\kappa}$, but also asserts that a significant portion of this variance occurs inside the event $|\a| \leq \kappa$, and also contains some more technical phase information about the covariance matrix of $\Re(\a)$ and $\Im(\a)$.

The following lemma was established in \cite{TVcir}:

\begin{lemma}\label{rot}  \cite[Lemma 2.4]{TVcir} Let $\a$ be a complex random variable with finite non-zero variance.  Then there exists a phase $e^{i\theta}$ and a $\kappa \geq 1$ such that $e^{i \theta} \a$ has $\kappa$-controlled second moment.
\end{lemma}

Since rotation by a phase does not affect the conclusion of Theorem \ref{theorem:main}, we conclude that we can assume without loss of generality that $\a$ is $\kappa$-controlled for some $\kappa$.  This will allow us to invoke several estimates from \cite{TVcir} (e.g. Lemma \ref{jbound} and Theorem \ref{theorem:TV} below).

\begin{remark} The estimates we obtain for Theorem \ref{theorem:main} will depend on $\kappa$ but will not otherwise depend on the precise distribution of $\a$.  It is in fact quite likely that the results in this paper can be generalised to random matrices $N_n$ whose entries are independent and are all $\kappa$-controlled for a single $\kappa$, but do not need to be identical.  In order to simplify the exposition, however, we focus on the iid case.
\end{remark}

\section{Small ball bounds} \label{section:smallball}

In this section we give some bounds on the small ball probabilities $\P(|\xi_1v_1 + \dots+ \xi_n v_n -z| \leq \eps)$ under various assumptions on the random variables $\xi_i$ and the coefficients $v_i$.  As a consequence we shall be able to obtain good bounds on the probability that $Av$ is small, where $A$ is a random matrix and $v$ is a fixed unit vector.

We first recall a standard bound (cf. \cite[Lemmas 4.2, 4.3, 5.2]{TVcir}):

\begin{lemma}[Fourier-analytic bound]\label{lemma:variance} Let $\xi_1, \dots, \xi_n$ be
independent variables.  Then we have the bound
$$\P(|\xi_1v_1 + \dots+ \xi_n v_n -z| \leq r) \ll r^2 \int_{w \in \C: |w| \leq 1/r} \exp( - \Theta( \sum_{j=1}^n \|w v_j\|_j^2  ) )\ dw $$
for any $r > 0$ and $z \in \C$, and any unit vector $v = (v_1,\ldots,v_n)$, where
\begin{equation}\label{jdef}
 \|z\|_j := ( \E \| \Re( z( \xi_j - \xi'_j) ) \|_{\R/\Z}^2)^{1/2},
 \end{equation}
$\xi'_j$ is an independent copy of $\xi_j$, and $\|x\|_{\R/\Z}$ denotes the distance from $x$ to the nearest integer.
\end{lemma}

\begin{proof}
By the Ess\'een concentration inequality (see e.g. \cite[Lemma 7.17]{TVbook}), we have
$$
\P(|\xi_1v_1 + \dots+ \xi_n v_n -z| \leq r )
\ll r^2 \int_{w \in \C: |w| \leq 1/r} |\E( e(\Re( w(\xi_1v_1 + \dots+ \xi_n v_n) ) ) )|\ dw$$
for any $c > 0$, where $e(x) := e^{2\pi i x}$.  We can write the right-hand side as
$$ r^2 \int_{w \in \C: |w| \leq 1/r} \prod_{j=1}^n f_j( w v_j )^{1/2}\ dw$$
where
$$ f_j(z) := |\E( e( \Re(\xi_j z) ) )|^2 = \E \cos(2\pi \Re( z(\xi_j - \xi'_j) ) ).$$
Using the elementary bound $\cos(2\pi \theta) \leq 1 - \Theta(\|\theta\|_{\R/\Z}^2)$ we conclude
$$ f_j(z) \leq 1 - \Theta( \| z \|_j^2 ) \leq \exp( - \Theta( \| z \|_j^2 ) )$$
and the claim follows.
\end{proof}

Next, we recall some properties of the norms $\| z\|_j$ in the case when $\xi_j$ is $\kappa$-controlled.

\begin{lemma}\label{jbound}  Let $1 \leq j \leq n$, let $\xi_j$ be a random variable, and let $\| \|_j$ be defined by \eqref{jdef}.
\begin{itemize}
\item[(i)] For any $w \in \C$, $0 \leq \|w\|_j \leq 1$ and $\|-w\|_j = \|w\|_j$.
\item[(ii)] For any $z,w \in \C$, $\|z+w\|_j \leq \|z\|_j+\|w\|_j$.
\item[(iii)] If $\xi_j$ is $\kappa$-controlled for some fixed $\kappa$,
then for any sufficiently small positive constants $c_0, c_1> 0$ we
have $\| z \|_j \ge c_1 \Re(z)$ whenever $|z| \leq c_0$.
\end{itemize}
\end{lemma}

\begin{proof} See \cite[Lemma 5.3]{TVcir}.
\end{proof}

We now use these bounds to estimate small ball probabilities.  We begin with a crude bound.

\begin{corollary}\label{cor:variance} Let $\xi_1, \dots, \xi_n$ be
independent variables which are $\kappa$-controlled.  Then there exists a constant $c > 0$ such that
\begin{equation}\label{Xiv}
\P(|\xi_1v_1 + \dots+ \xi_n v_n -z| \leq c) \leq 1-c
\end{equation}
for all $z \in \C$ and all unit vectors $(v_1,\ldots,v_n)$.
\end{corollary}

\begin{proof} Let $c > 0$ be a small number to be chosen later.  We divide into two cases,
depending on whether all the $v_i$ are bounded in magnitude by $\sqrt{c}$ or not.

Suppose first that $|v_i| \leq \sqrt{c}$ for all $c$.  Then we apply
 Lemma \ref{lemma:variance} (with $r := c^{1/4}$) and bound the left-hand side of \eqref{Xiv} by
$$ \ll c^{1/2} \int_{w \in \C: |w| \leq c^{-1/4}} \exp( - \Theta( \sum_{j=1}^n \|w v_j\|_j^2  ) )\ dw.$$
By Lemma \ref{jbound}, if $c$ is sufficiently small then we have $\|
w v_j \|_j \ge c_1 \Re(w v_j)$, for some positive constant $c_1$.
Writing each $v_j$ in polar coordinates as $v_j = r_j e^{2\pi i
\theta_j}$, we thus obtain an upper bound of
$$ \ll c^{1/2} \int_{w \in \C: |w| \leq c^{-1/4}} \exp( - \Theta( \sum_{j=1}^n r_j^2 \Re(e^{2\pi i \theta_j} w)^2  ) )\ dw.$$
Since $\sum_{j=1}^n r_j^2 = 1$, we can use H\"older's inequality (or Jensen's inequality) and bound this from above by
$$ \ll \sup_j c^{1/2} \int_{w \in \C: |w| \leq c^{-1/4}} \exp( - \Theta( \Re(e^{2\pi i \theta_j} w)^2  ) )\ dw$$
which by rotation invariance and scaling is equal to
$$ \int_{w \in \C: |w| \leq 1} \exp( - \Theta( c^{-1/4} \Re(w)^2  ) )\ dw.$$
From the monotone convergence theorem (or direct computation) we see
that this quantity is less than $1-c$ if $c$ is chosen sufficiently
small. (If necessary, we allow $c$ to depend on the hidden constant
in $\Theta$.)

Now suppose instead that $|v_1| > \sqrt{c}$ (say).
 Then by freezing all of the variables $\xi_2,\ldots,\xi_n$, we can bound the left-hand side of \eqref{Xiv} by
$$ \sup_w \P( |\xi_1 - w| \leq \sqrt{c} ).$$
But by the definition of $\kappa$-control,  one easily sees that
this quantity is bounded by $1-c$ if $c$ is sufficiently small
(compared to $1/\kappa$), and the claim follows.
\end{proof}

As a consequence of this bound, we obtain

\begin{theorem} \label{theorem:largest3} Let $N_n$ be an $n \times
n$ random matrix whose entries are independent random variables which are all
 $\kappa$-controlled for some
 constant $\kappa > 0$.
  Then there are positive constants $c, c'$ such that the following holds.
   For any unit vector $v$ and any deterministic matrix
  $M$,  $$\P(\| (M+N_n)v \| \le cn^{1/2} ) \le \exp(-c'n).$$
\end{theorem}

\begin{proof}
Let $c$ be a sufficiently small constant, and let $X_1,\ldots,X_n$
denote the rows of $M+N_n$. If $\|(M+N_n)v \| \le c n^{{1/2}}$, then
we have $|\langle X_j,v\rangle | \le c$ for at least $(1-c)n$ rows.
As the events $\I_j:= |\langle X_j,v\rangle | \le c$ are
independent, we see from the Chernoff inequality (applied to the sum
$\sum_j \I_j$ of indicator variables)  that it suffices to show that
$$ \E(I_j)= \P( |\langle X_j,v\rangle | \le c ) \leq 1-2c$$
(say) for all $j$.  But this follows from Corollary \ref{cor:variance} (after adjusting $c$ slightly), noting that each $X_j$ is a translate (by a row of $M$) of a vector whose entries are iid copies of $\a$.
\end{proof}

Now we obtain some statements
 of inverse Littlewood-Offord type.

\begin{definition}[Compressible and incompressible vectors]  For any $a, b > 0$,
let $\Comp(a,b)$ be the set of unit vectors $v$ such that there is a
vector $v'$ with at most $an$ non-zero coordinates satisfying $\| v-v' \| \le b$.
We denote by $\Incomp (a,b)$ the set of unit vectors which do not lie in $\Comp(a,b)$.
\end{definition}

\begin{definition}[Rich vectors]\label{rich}  For any $\eps, \rho > 0$, let $S_{\eps, \rho}$ be the set of unit vectors $v$  satisfying
$$ \sup_{z \in \C}  \P ( |X \cdot v -z | \le \eps)  \ge \rho,$$
where $X = (\a_1,\ldots,\a_n)$ is a vector whose coefficients are iid copies of $\a$.
\end{definition}

\begin{lemma}[Very rich vectors are compressible]\label{lemma:lastcase}
For any $\eps, \rho > 0$ we have
$$ S_{\eps, \rho} \subset \Comp\left( O( \frac{1}{n\rho^2} ), O( \frac{\eps}{\rho} ) \right).$$
\end{lemma}

\begin{proof}  We can assume $\rho \gg n^{-1/2}$ since the claim is trivial otherwise.
Let $v \in S_{\eps, \rho}$, thus
$$ \P(|X \cdot v-z| \le \eps ) \ge \rho $$
for some $z$.
From Lemma \ref{lemma:variance} we conclude
\begin{equation}\label{nancy}
 \eps^2 \int_{w \in \C: |w| \leq \eps^{-1}} \exp( - \Theta( \sum_{j=1}^n \|w v_j\|_j^2  ) )\ dw \gg \rho.
\end{equation}

Let $s > 0$ be a small constant (independent of $n$) to be chosen later, and let $A$ denote the set of indices $i$ for which $|v_i| \geq s \eps$.  Then from \eqref{nancy} we have
$$ \eps^2 \int_{w \in \C: |w| \leq \eps^{-1}} \exp( - \Theta( \sum_{j \in A} \|w v_j\|_j^2  ) )\ dw \gg \rho.$$
Suppose $A$ is non-empty.  Applying H\"older's inequality, we conclude that
$$  \eps^2 \int_{w \in \C: |w| \leq \eps^{-1}} \exp( - \Theta( |A| \|w v_j\|_j^2 ) )\ dw \gg \rho $$
for some $j \in A$. By the pigeonhole principle, this implies that
\begin{equation}\label{pigeon}
| \{ w \in \C: |w| \leq \eps^{-1}, |A| \|w v_j \|_j^2 \leq k \} | \gg k^{1/2} \eps^{-2} \rho
\end{equation}
for some integer $k \geq 1$.

If $|A| \ll k$, then the set in \eqref{pigeon} has measure
$\Theta(\eps^{-2})$, which forces $|A| \ll \rho^{-2}$.  Suppose
instead that $k \leq s |A|$ for some small $s' > 0$. Since
$| v_j| \ge s \epsilon$, we have $s'/|v_j| \le
s'/s \epsilon$. We will choose $s'$ sufficiently
small to make sure that this ratio is smaller than the constant
$c_0$ in Lemma \ref{jbound}. By Lemma \ref{jbound}, we see that the
intersection of the set in \eqref{pigeon} with any ball of radius
$s'/|v_j|$ has density at most $\sqrt{k/|A|}$, and so by
covering arguments we can bound the left-hand side of \eqref{pigeon}
from above by $\ll k^{1/2} |A|^{-1/2} \eps^{-2}$. Thus we have $|A|
\ll \rho^{-2}$ in this case also.  Thus we have shown in fact that
$|A| \ll \rho^{-2}$ in all cases (the case when $A$ is empty being
trivial).

Now we consider the contribution of those $j$ outside of $A$.  From \eqref{nancy} and Lemma \ref{jbound} we have
$$
 \eps^2 \int_{w \in \C: |w| \leq \eps^{-1}} \exp( - \Theta( \sum_{j \not \in A} \Re( w v_j )^2  ) )\ dw \gg \rho.$$
Suppose that $A$ is not all of $\{1,\ldots,n\}$.  Using polar coordinates $v_j = r_j e^{2\pi i \theta_j}$ as before, we see from H\"older's inequality that
$$
 \eps^2 \int_{w \in \C: |w| \leq \eps^{-1}} \exp( - \Theta( r^2 \Re( w e^{2\pi i \theta_j} )^2  ) )\ dw \gg \rho
 $$
for some $j \not \in A$, where $r^2 := \sum_{j \not \in A} r_j^2$.  After scaling and rotation invariance, we conclude
$$
 \int_{w \in \C: |w| \leq 1} \exp( - \Theta( \frac{r^2}{\eps^2} \Re( w )^2  ) )\ dw \gg \rho.
$$
The left-hand side can be computed to be at most $O(\eps/r)$.  We conclude that $r \ll \eps/\rho$.  If we let $v'$ be the restriction of $v$ to $A$, we thus have $\|v-v'\| \ll \eps/\rho$, and the claim $v \in
\Comp( O( \frac{1}{n\rho^2} ), O( \frac{\eps}{\rho} ) )$ follows.  (The case when $A=\{1,\ldots,n\}$ is of course trivial.)
\end{proof}

Roughly speaking, Lemma \ref{lemma:lastcase} gives a complete characterization of
vectors $v$ such that $$ \sup_{z \in \C}  \P ( |X \cdot v -z | \le
\eps) \ge \rho,$$ where $\rho > C n^{-1/2} $, for some large
constant $C$. The lemma shows that such a vector $v$ can be
approximated by a vector $v'$  with at most $\frac{C'}{\rho^2}$
non-zero coordinates such that $\|v-v'\| \le \frac{C^{''}
\epsilon}{\rho}$, where $C', C^{''}$ are positive constants.

The dependence of parameters here are sharp, up to constant terms. Indeed, in the
Bernoulli case, the vector $v = (1,\ldots,1,0,\ldots,0)$ consisting
of $k$ $1$s lies in $S_{0,\Theta(1/\sqrt{k})}$ and lies in $\Comp(
a, 0)$ precisely when $an \geq k$ (cf. \cite{erdos}). This shows
that the $O(\frac{1}{n\rho^2})$ term on the right-hand side cannot
be improved.  On the other hand, in the Gaussian case, observe that
if $\|v\| \leq b$ then $X \cdot v$ will have magnitude $O(\eps)$
with probability $O( \eps/b )$, which shows that the term
$O(\frac{\eps}{\rho})$ cannot be improved.

Lemma \ref{lemma:lastcase} is only non-trivial in the case $\rho \ge
C n^{-1/2}$, for some large constant $C$.  To handle the case of
smaller $\rho$, we use the following more difficult entropy bound
from \cite{TVcir}.

\begin{theorem}[Entropy of rich vectors]\label{theorem:TV}
 For any $\eps, \rho$, there is a finite set $S'_{\eps, \rho}$ of size at most
 $n^{-(1/2-o(1)) n} \rho^{-n}  + \exp(o(n))$ such that for each
 $v \in S_{\eps, \rho}$, there is
 $v' \in S'_{\eps, \rho}$ such that $ \| v-v' \|_{\infty} \le \eps$.
\end{theorem}

\begin{proof} See \cite[Theorem 3.2]{TVcir}. \end{proof}

\section{Proof of Theorem \ref{theorem:main}: preliminary reductions}

We now begin the proof of Theorem \ref{theorem:main}.
Let $N_n, M, \gamma, A$ be as in that theorem. As remarked in Section \ref{section:generaltheorems},
we may assume $\a$ to be $\kappa$-controlled for some $\kappa$.
We allow all implied constants to depend on $\kappa, \gamma, A$.
We may of course assume that $n$ is large compared to these parameters.  We may also assume that
\begin{equation}\label{n2}
\P( \|N_n\| \geq n^\gamma ) \leq \frac{1}{2}
\end{equation}
since the claim is trivial otherwise. By decreasing $A$ if necessary, we may furthermore assume that
\begin{equation}\label{na}
\P( \|N_n\| \geq n^\gamma ) \leq n^{-A+o(1)}.
\end{equation}

It will then suffice to show (assuming \eqref{n2}, \eqref{na}) that
$$\P (s_n(M_n) \leq n^{-(2A+1)\gamma}) \ll n^{-A+\alpha+o(1)}$$
for any constant  $\alpha > 0$ (with the implied constants now
depending on $\alpha$ also), since the claim then follows by sending
$\alpha$ to zero very slowly in $n$.

Fix $\alpha$, and allow all implied constants to depend on $\alpha$.
By perturbing $A$ and $\alpha$ slightly we may assume that $A$ is
not a half-integer; we can also take $\alpha$ to be small depending
on $A$. For example, we can assume that

\begin{equation} \label{eqn:alphaA} \alpha < \{2A \} /2
\end{equation}

\noindent where $\{2A\}$ is the fractional part of $2A$.

Using the trivial bound $\| N_n \| \geq \sup_{1 \leq i,j \leq n} |\a_{ij}|$, we conclude from \eqref{n2}, \eqref{na} that
$$
\P( |\a_{ij}| \geq n^\gamma \hbox{ for some } i,j ) \leq \min(\frac{1}{2}, n^{-A+o(1)}).
$$
Since $\a_{ij}$ are iid copies of $\a$, the  $n^2$ events $|\a_{ij}|
\geq n^\gamma$ are independent with identical probability. It
follows that
\begin{equation}\label{ngamma}
\P( |\a| \geq n^\gamma ) \leq n^{-A-2+o(1)}.
\end{equation}

Let $F$ be the event that $s_n(M_n) \leq n^{-(2A+1)\gamma}$, and
let $G$ be the event that $\|N_n\| \leq n^\gamma$.  In view of \eqref{na}, it suffices to show that
$$\P( F \wedge G ) \le n^{-A+\alpha+o(1)}.$$
Set
\begin{equation}\label{bdef}
b:= \beta n^{1/2-\gamma}
\end{equation}
and
\begin{equation}\label{Adef}
a := \frac{\beta}{\log n },
\end{equation}
where $\beta$ is a small positive constant  to be chosen later. We
then introduce the following events:

\begin{itemize}
\item $F_{\Comp}$ is the event that $\|M_n v \| \leq n^{-(2A+1)\gamma}$ for some $v \in \Comp(a,b)$.
\item $F_{\Incomp}$ is the event that $\|M_n v \| \leq n^{-(2A+1)\gamma}$ for some $v \in \Incomp(a,b)$.
\end{itemize}

Observe that if $F$ holds, then at least one of $F_{\Comp}$ and
$F_{\Incomp}$ holds. Theorem \ref{theorem:main} then follows
immediately from the following two lemmas.

\begin{lemma}[Compressible vector bound]\label{lemma:comp} If $\beta$ is sufficiently small,
then $$\P(F_{\Comp} \wedge G) \le \exp(-\Omega(n)). $$
\end{lemma}

 \begin{lemma}[Incompressible vector bound]\label{lemma:incomp} We have
 $$\P(F_{\Incomp} \wedge G) \le n^{-A+o(1)}. $$
 \end{lemma}

In these lemmas we allow the implied constants to depend on $\beta$.

The proof of Lemma \ref{lemma:comp} is simple and will be presented in the next section. The proof of Lemma \ref{lemma:incomp} is somewhat more involved and occupies the rest of the paper.

\section{Treatment of compressible vectors}
 \label{section:proofcomp}

If $F_{\Comp} \wedge G$ occurs, then by the definition of
$\Comp(a,b)$,
 there are unit vectors $v, v'$ such that $\| M_n v\| \le n^{-(2A+1)\gamma}$ and $v'$ has support on at
 most  $an$ coordinates and $\|v-v'\| \le b$.

  By the triangle inequality and \eqref{bdef} we have
\begin{align*}
\| M_nv' \| &\le n^{-(2A+1)\gamma} +  \|M_n\|  \|v-v' \| \\
&\le n^{-(2A+1)\gamma} + n^{\gamma}b\\
&\leq 2 \beta n^{1/2}.
\end{align*}

A set $\N$ of unit vectors in $\C^m$ is called a \emph{$\delta$-net} if for
any unit vector $v$, there is a vector $w$ in $\N$ such that
$\|v-w\| \le \delta$. It is well known that for any $0 < \delta <1$,
a $\delta$-net of size $(C\delta^{-1})^m$ exists, for some constant
$C$ independent of $\delta$ and $m$.

Using this fact, we conclude that the set  of unit vectors with  at
most $an$ non-zero coordinates admits an $b$-net $\N$ of size at
most
$$ |\N| \le \binom{n}{an} (C b^{-1} )^{an},$$

\noindent  Thus, if $F_{\Comp} \wedge G$ occurs, then there is a
unit vector $v'' \in \N$   such that

$$\|M_n  v'' \| \le 2 \beta n^{1/2}  + \|M_n \| b =  3 \beta n^{1/2}. $$

On the other hand, from Theorem \ref{theorem:largest3} we see (for
$\beta \le c/3$) that for any fixed $v''$, $$\P( \| M_n v'' \| \leq
3 \beta n^{1/2} ) \leq \exp(-c'n), $$  where $c$ and $c'$ are the
constants in Theorem \ref{theorem:largest3}.

 By the union bound, we conclude
$$ \P(F_{\Comp} \wedge G) \le \binom{n}{an}  (b^{-1} )^{an} \exp(-c'n).$$

But from \eqref{bdef}, \eqref{Adef} we see that the right-hand side
can be made less than  $\exp(-c'n/2)$, given that  $\beta$ is
sufficiently small. This concludes the proof of Lemma
\ref{lemma:comp}.

\section{Treatment of incompressible vectors} \label{section:proofincomp}

We now begin the proof of Lemma \ref{lemma:incomp}.
 We now fix $\beta$ and allow all implied constants to depend on $\beta$.

Let $X_k$ be the $k^{\th}$ row vector of $M_n$, and let $\dist_k$ be
the distance from $X_k$ to the subspace spanned by $X_1, \dots,
X_{k-1}, X_{k+1}, \dots, X_n$.  We need the following, which is a
slight extension of a  lemma from \cite{RV}.

\begin{lemma} \label{lemma:RV} For any $\eps > 0$, and any event $E$, we have
$$ \P (\{\|Mv\| \le \eps b n^{-1/2} \hbox{ for some } v \in \Incomp(a,b) \} \wedge E
) \le \frac{1}{an} \sum_{k=1}^{n} \P (\{ \dist_k \le \eps \}\wedge
E).
$$
\end{lemma}

\begin{proof} See \cite[Lemma 3.5]{RV}.  The arbitrary event $E$
was not present in that lemma, but one easily verifies that the proof works perfectly well with this event in place.
\end{proof}

Applying this to our current situation with
\begin{equation}\label{epsdef}
\eps := \frac{1}{\beta} n^{-2A\gamma},
\end{equation}
we obtain
$$ \P( F_{\Incomp} \wedge G ) \ll \frac{\log n}{n} \sum_{k=1}^n \P( \{\dist_k \leq \eps\} \wedge G ).$$
To prove Lemma \ref{lemma:incomp}, it therefore suffices (by symmetry) to show that
$$ \P( \{\dist_n \leq \eps \} \wedge G) \ll n^{-A+\alpha+o(1)}.$$
Notice that there  is a unit vector $X_n^{\ast}$ orthogonal to
$X_1, \dots, X_{n-1}$ such that
\begin{equation}\label{disk}
\dist_k = |X_n\cdot X_n^{\ast} | .
\end{equation}
\noindent If there are many such $X_n^{\ast}$, choose one
arbitrarily.  However, note that we can choose $X_n^\ast$ to depend only on $X_1,\ldots,X_{n-1}$ and thus be independent of $X_n$.

Let $\rho := n^{-A+\alpha}$. Let $X$ be the random vector of length
$n$ whose coordinates are iid copies of $\a$.
 From Definition \ref{rich} (and the observation that $X_n$ has the same distribution
 as $X$ after translating by a deterministic vector (namely the $n$th row of the deterministic matrix $M$),
  we have the conditional probability bound
$$ \P( \dist_n \leq \eps | X_n^\ast \not \in S_{\eps,\rho}) \leq \rho = n^{-A+\alpha}.$$
Thus it will suffice to establish the exponential bound

$$ \P( \{X_n^\ast \in S_{\eps,\rho} \} \wedge G ) \le \exp( - \Omega(n) ).$$

Let
\begin{equation}\label{J-def}
J := \lfloor 2A \rfloor
\end{equation}
be the integer part of $2A$. Let $\alpha_1>0$ be a sufficiently
small constant (independent of $n$ and $\gamma$, but depending on
$\alpha, A, J$) to be chosen later.  Set
\begin{equation}\label{epsj-def}
 \eps_j := n^{(\gamma+\alpha_1)j}\eps = \frac{1}{\beta} n^{(\gamma+\alpha_1)j} n^{-2A\gamma}
 \end{equation}
and
\begin{equation}\label{rhoj-def}
\rho_j := n^{(1/2-\alpha_1)j} \rho = n^{(1/2-\alpha_1)j} n^{-A+\alpha}
\end{equation}
for all $0 \leq j \leq J$.

By the union bound, it will suffice to prove the following lemmas.

\begin{lemma}\label{j1} If $\alpha_1$ is sufficiently small, then for any $0 \leq j < J$, we have
\begin{equation}\label{pxn}
\P( \{X_n^\ast \in S_{\eps_j, \rho_j}\} \wedge \{X_n^\ast \not \in
S_{\eps_{j+1}, \rho_{j+1}}\} \wedge G ) \le \exp( - \Omega(n) ).
\end{equation}
\end{lemma}

\begin{lemma}\label{j2} If $\alpha_1$ is sufficiently small, then we have
$$ \P( X_n^\ast \in S_{\eps_J, \rho_J} ) \le \exp( - \Omega(n) ).$$
\end{lemma}

\section{Proof of Lemma \ref{j1}} \label{prooflemma:j1}

Fix $0 \leq j < J$.  Note that by \eqref{eqn:alphaA}, we have

$$ \rho_j \leq n^{(J-1)/2} n^{-A+\alpha} \leq n^{-1/2 - \{2A\}/2 + \alpha} \leq n^{-1/2}.$$

 We can then use Theorem \ref{theorem:TV} to conclude the existence of
a set $\N$ of unit vectors such that every vector in $S_{\eps_j,
\rho_j}$ lies within $\eps_j$ in $l^\infty$ norm to a vector in
$\N$, and with the cardinality bound
\begin{equation}\label{noon}
 |\N| \le n^{-(1/2-o(1)) n} \rho_j^{-n}.
 \end{equation}

Suppose that the event in Lemma \ref{j1} holds, then we can find $u
\in \N$ such that $\| u-X_n^\ast\|_{l^\infty} \leq \eps_j$, and thus
$\| u - X_n^\ast \| \leq n^{1/2} \eps_j$. On the other hand, since
$X_n^\ast$ is orthogonal to $X_1,\ldots,X_{n-1}$ and $\|M_n\| \ll
n^\gamma$, we have
\begin{align*}
(\sum_{i=1}^{n-1} |X_i \cdot u|^2)^{1/2} &= (\sum_{i=1}^{n-1} |X_i
\cdot (u-X_n^\ast)|^2)^{1/2} \\
&= \| M (u-X_n^\ast)\| \\
&\ll n^\gamma n^{1/2} \eps_j\\
&\ll n^{1/2} n^{-\alpha_1} \eps_{j+1}.
\end{align*}
On the other hand, from \eqref{pxn} and Definition \ref{rich} we have
\begin{equation}\label{xxn}
 \P( |X \cdot X_n^\ast - z| \leq \eps_{j+1} ) \leq \rho_{j+1}
\end{equation}
for all $z \in \C$, where $X = (\a_1,\ldots,\a_n)$ consists of iid copies of $\a$.

To conclude the proof, we will need the following lemma.

\begin{lemma} \label{lemma:momentmethod} If $w$ is any vector with $\|w\|_{l^\infty} \leq 1$, then
$$ \P( |X \cdot w| \geq n^{\gamma+\alpha_1} ) \ll n^{-A}.$$
\end{lemma}

\begin{proof}  Write $w = (w_1,\ldots,w_n)$ and $X = (\a_1,\ldots,\a_n)$.
Observe from \eqref{na} that with probability $O(n^{-A-1}) = O(n^{-A})$,
all the coefficients in $X$ are going to be of magnitude at most $n^\gamma$.  Thus it suffices to show that
$$ \P( |w_1 \tilde \a_1 + \ldots + w_n \tilde \a_n| \geq n^{\gamma + \alpha_1} ) \ll n^{-A} $$
where $\tilde \a_1,\ldots,\tilde \a_n$ are iid with law equal to that of $\a$ conditioned to the event $|\a| \ll n^\gamma$.  As $\a$ has mean zero and bounded second moment, one verifies from \eqref{na} and Cauchy-Schwarz that the mean of the $\tilde \a_i$ is $O(n^{-(A+2)/2})$.  Thus if we let $\a'_i := \tilde \a_i - \E( \tilde \a_i)$, we see that it suffices to show that
$$ \P( |w_1 \a'_1 + \ldots + w_n \a'_n| \geq \frac{1}{2} n^{\gamma + \alpha_1} ) \ll n^{-A}.$$

We conclude the proof by the moment method, using the following
estimate
$$ \E( |w_1 \a'_1 + \ldots + w_n \a'_n|^{2k} ) \ll_k n^{2k \gamma}$$
for any integer $k \geq 0$. This is easily verified by a standard
computation (using the hypothesis $\gamma \geq 1/2$), since all the
$\a'_i$ have vanishing first moment, a second moment of $O(1)$, and
a $j^{th}$ moment of $O_j( n^{(j-2) \gamma} )$ for any $j > 2$. Now
take $k$ to be a constant sufficiently large compared to
$A/\alpha_1$.
\end{proof}

We are now ready to finish the proof of Lemma \ref{j1}.  From lemma
\ref{lemma:momentmethod}  and the bound $\|u-X_n^\ast\| \leq \eps_j$
we see that
$$ \P( |X \cdot (X_n^\ast-u)| \geq \eps_{j+1} ) \leq n^{-A} \leq \rho_{j+1};$$
combining this with \eqref{xxn} using the triangle inequality, we see that
\begin{equation}\label{xuz}
\sup_{z \in \C} \P( |X \cdot u - z| \leq \eps_{j+1} ) \ll \rho_{j+1}.
\end{equation}
We can therefore bound the left-hand side of \eqref{pxn} by
$$ \sum_{u \in \N: \eqref{xuz} \hbox{ holds}}
\P \Big( (\sum_{i=1}^{n-1} |X_i \cdot u|^2)^{1/2} \ll n^{1/2}
n^{-\alpha_1} \eps_{j+1} \Big).$$ Now suppose that $u \in \N$ obeys
\eqref{xuz}.  If we have $\sum_{i=1}^{n-1} |X_i \cdot u|^2)^{1/2}
\ll n^{1/2} n^{-\alpha_1} \eps_{j+1}$, then the event $|X_i \cdot u|
\leq \eps_{j+1}$ must hold for at least $n - O( n^{1-2\alpha_1})$
values of $i$.  On the other hand, from \eqref{xuz} we see that each
of these events $|X_i \cdot u| \leq \eps_{j+1}$ only occurs with
probability $O(\rho_{j+1})$.  We can thus bound
\begin{align*}
\P( \sum_{i=1}^{n-1} |X_i \cdot u|^2)^{1/2} \ll n^{1/2} n^{-\alpha_1} \eps_{j+1} )
&\leq \binom{n}{n - O( n^{1-2\alpha_1})} ( O(\rho_{j+1}) )^{n - O( n^{1-2\alpha_1})} \\
&\ll n^{o(n)} \rho_{j+1}^n.
\end{align*}
Applying \eqref{noon}, we can thus bound the left-hand side of \eqref{pxn} by
$$ \ll n^{-(1/2-o(1)) n} \rho_j^{-n} \rho_{j+1}^n = n^{-(\alpha_1-o(1))n}$$
and the claim follows.

\section{Proof of Lemma \ref{j2}} \label{prooflemma:j2}

Suppose that $X_n^\ast$ lies in $S_{\eps_J, \rho_J}$.  Then by Lemma \ref{lemma:lastcase}, we have
$$ X_n^\ast \subset \Comp( O( \frac{1}{n\rho_J^2} ), O( \frac{\eps_J}{\rho_J} ) ).$$
Note from \eqref{rhoj-def} and \eqref{J-def} that
$$ \frac{1}{n\rho_J^2} = n^{2A-J-1+2\alpha_1 J-2\alpha} \leq n^{-\alpha_1}$$
if $\alpha_1$ is sufficiently small.  Thus, by arguing as in Section \ref{section:proofcomp}, the set $\Comp( O( \frac{1}{n\rho_J^2} ), O( \frac{\eps_J}{\rho_J} ) )$ has a $O( \frac{\eps_J}{\rho_J} )$-net $\N$ in $l^2$ of cardinality
$$ |\N| \ll \binom{n}{\frac{1}{n\rho_J^2}} (O( \frac{\eps_J}{\rho_J} ))^{\frac{1}{n\rho_J^2}} = \exp(o(n)).$$
If we let $u \in \N$ be within $O( \frac{\eps_J}{\rho_J} )$ of $X_n^\ast$, then we have $|X_i \cdot u| \ll \frac{\eps_J}{\rho_J}$ for all $1 \leq i \leq n-1$.  Thus we can bound
$$
\P( X_n^\ast \in S_{\eps_J, \rho_J} ) \leq \sum_{u \in \N} \P( |X_i \cdot u| \ll \frac{\eps_J}{\rho_J} \hbox{ for all } 1 \leq i \leq n-1 ).$$
Now observe from \eqref{epsj-def}, \eqref{rhoj-def}, \eqref{J-def} and the hypothesis $\gamma \geq 1/2$ that
$$ \frac{\eps_J}{\rho_J} = n^{-\alpha+2\alpha_1 J} n^{-(2A-J)(\gamma-1/2)} \leq n^{-\alpha/2}$$
(say) if $\alpha_1$ is sufficiently small.  Thus by Corollary \ref{cor:variance} (or by a minor modification of Theorem \ref{theorem:largest3}) we see that
$$ \P( |X_i \cdot u| \ll \frac{\eps_J}{\rho_J} \hbox{ for all } 1 \leq i \leq n-1 ) \ll \exp(-\Omega(n))$$
for each $u \in \N$, and the claim follows.

\vskip2mm

{\it Acknowledgement.} We would like to thank the referees for useful comments.


\begin{thebibliography}{10}

\bibitem{BS}  Z. Bai and J. Silverstein, Spectral analysis of large dimensional random matrices, Science Press, 2007.

\bibitem{BT} D. Bau and L. N. Trefethen, Numerical linear algebra,
{SIAM 1997}.


\bibitem{BCL0} P. Burgisser, F. Cucker, M. Lotz, The probability that a slightly perturbed numberical
analysis problem is difficult, {\it Math of Computation,} 77,
1559-1583, 2008.


\bibitem{BCL1} P. Burgisser, F. Cucker, M. Lotz, General formulas for the smooth analysis of condition numbers,
{\it C. R. Acad. Sc. Paris}, 343, 145-150, 2006.

\bibitem{BCL2} P. Burgisser, F. Cucker, M. Lotz, Smooth analysis of connic condition numbers,
{\it J. Math. Pure et Appl.}, 86, 293-309, 2006.




\bibitem{Edel}  A. Edelman,  Eigenvalues and condition numbers of random matrices.
{\it SIAM J. Matrix Anal. Appl.}  9 (1988), no. 4, 543--560.

\bibitem{erdos}
P. Erd\"os, On a lemma of Littlewood and Offord, {\it  Bull. Amer.
Math. Soc.}  \textbf{51} (1945), 898--902.

\bibitem{Demmel} J. Demmel, The probability that a  numberical analysis problem is difficult, {\it Math. Comp.}, 50, 449-480, 1988

\bibitem{DST} J. Dunagan, D. A. Spielman  and S. H. Teng,
Smoothed Analysis of the Renegar's Condition Number for Linear
Programming, {\it preprint}.

\bibitem{GV} G. Golub and C. Van Loan,  Matrix computations, Third edition, Johns
Hopkins Press 1996.

\bibitem{GvN}
J. von Neumann, H. Goldstine, \emph{Numerical inverting of matrices of high order},
Bull. Amer. Math. Soc. \textbf{53} (1947). 1021--1099.

\bibitem{FK} Z. F\"uredi and J. Koml\'os, The eigenvalues of random symmetric matrices,{\it Combinatorica } 1 (1981), no. 3, 233--241.



\bibitem{Kos} Kostlan,

\bibitem{Kom} J. Koml\'os, On the determinant of $(0,1)$ matrices,
  {\it Studia Sci. Math. Hungar.} \textbf{ 2}
(1967) 7-22.

\bibitem{Latala} R. Latala, Some estimates of norms of random matrices, \emph{Proc. Amer. Math. Soc.} \textbf{133} (2005), 1273--1282.

\bibitem{Lit} A. Litvak, A. Pajor, M. Rudelson and N. Tomczak-Jaegermann,
 Smallest singular value of random matrices and geometry of random polytopes,
 {\it  Adv. Math. } 195 (2005), no. 2, 491--523.

\bibitem{Mehta} M.L. Mehta, Random Matrices and the Statistical Theory of Energy Levels, Academic Press, New York, NY, 1967.




\bibitem{Pas} L. A Pastur, On the spectrum of random matrices, {\it Teoret. Mat. Fiz. } 10, 102-112 (1973).

\bibitem{Ren} J. Renegar, On the efficiency of Newton's method in approximating all zeros of a
system of complex polynomials, {\it Math. Oper. Res}, 12 (1987), no 1, 121-148.


\bibitem{Rud} M. Rudelson,  Invertibility of random matrices: Norm of the inverse.
     {\it Annals of Mathematics, to appear}.



\bibitem{RV} M. Rudelson and R. Vershynin, The Littlewood-Offord
problem and the condition number of random matrices, {\it Adv. Math., to appear}.

\bibitem{Smale}  S. Smale, On the efficiency of algorithms of analysis, {\it Bullentin of
the AMS} (13) (1985), 87-121.

\bibitem{Seg}  Y. Seginer, The expected norm of random matrices,
{\it   Combin. Probab. Comput. } 9  (2000),  no. 2, 149--166.



\bibitem{ST}  D. A. Spielman and S. H. Teng,
 Smoothed analysis of algorithms, {\it  Proceedings of the International Congress of
Mathematicians}, Vol. I (Beijing, 2002), 597--606, Higher Ed.
Press, Beijing, 2002.

\bibitem{ST1}  D. A. Spielman and S. H. Teng, Smoothed analysis of algorithms: why the simplex algorithm
usually takes polynomial time, {\it J. ACM } 51 (2004), no. 3,
385--463.

\bibitem{SST}  A. Sankar, S. H. Teng, and D. A. Spielman,
Smoothed Analysis of the Condition Numbers and Growth Factors of
Matrices,  {\it SIAM J. Matrix Anal. Appl.}  28  (2006),  no. 2, 446--476.




\bibitem{TVsing} T. Tao and V. Vu, Inverse Littlewood-Offord theorems and the condition number of
 random discrete matrices, {\it Annals of Mathematics, to appear}.

\bibitem{TVstoc} T. Tao and V. Vu, The condition number of a randomly perturbed matrix, {\it STOC 2007}.

\bibitem{TVcir} T. Tao and V. Vu, Random matrices: The circular law, {\it Communications in Contemporary Mathematics}, 10 (2008), 261-307.

\bibitem{TVbook} T. Tao and V.  Vu, Additive Combinatorics, Cambridge Univ. Press, 2006.

\bibitem{vonG} J. von Neuman and H. Goldstein, Numerical inverting matrices of high order,
{\it Bull. Amer. Math. Soc.} 53, 1021-1099, 1947.



\bibitem{Vnorm} V. Vu,
Spectral norm of random matrices, {\it  Combinatorica  27}   (2007),  no. 6, 721--736.


\bibitem{Wig} P. Wigner, On the distribution of the roots of certain symmetric matrices, {\it Annals of Math.}, 67, 325-327.


 \end{thebibliography}
 \end{document}


\title[Erratum to ``Smooth analysis...'']{Erratum to ``Smooth analysis of the condition number and the least singular value''}

\author{Terence Tao}
\address{Department of Mathematics, UCLA, Los Angeles CA 90095-1555}
\email{tao@math.ucla.edu}
\thanks{T. Tao is supported by NSF grant DMS-0649473 and by a Simons Investigator Award.}

\author{Van Vu}
\thanks{V. Vu is supported by research grants DMS-0901216 and AFOSAR-FA-9550-09-1-0167.}
\address{Department of Mathematics, Yale University, New Haven, CT 06520, USA}
\email{van.vu@yale.edu}

\begin{abstract}  We point out a gap in our previous paper ``Smooth analysis of the condition number and the least singular value'', which can be fixed at the cost of worsening the bound in the main theorem of that paper by a factor of $n^{\gamma-1/2}$.
\end{abstract}

\maketitle

\section{Introduction}

The purpose of this erratum is to point out a gap in our previous paper \cite{tv-cond}, which requires a modification to the statement of Theorem 3.2 of that paper (and hence also Corollaries 3.3, 3.5, 3.7, and 3.8).  We thank Nick Cook for alerting us to this issue.  The problem occurs in Section 11, in the line ``If we let $u \in {\mathcal N}$ be within $O(\frac{\eps_J}{\rho_J})$ of $X_n^*$, then we have $|X_i \cdot u| \ll \frac{\eps_J}{\rho_J}$''.  This implication is unjustified, and we regret this oversight.

We are able to repair the argument, but only at the cost of worsening the conclusion of \cite[Theorem 3.2]{tv-cond} from
$$ {\bf P}( s_n(M_n) \leq n^{-(2A+1)\gamma} ) \leq c \left( n^{-A+o(1)} + {\bf P}( \|N_n\| \geq n^\gamma ) \right)$$
to
$${\bf P}( s_n(M_n) \leq n^{-(2A+2)\gamma+1/2} ) \leq c \left( n^{-A+o(1)} + {\bf P}( \|N_n\| \geq n^\gamma ) \right),$$
thus the lower bound on $s_n(M_n)$ is worse by a factor of $n^{-(\gamma-1/2)}$.  This does not affect the most important case $\gamma= 1/2$, but causes some degradation in the numerology for $\gamma > 1/2$.  However, in all the applications of \cite{tv-cond} that we are aware of in the literature, the precise value of this exponent in the $\gamma>1/2$ case is not important to the analysis.

The repair to the proof of this modification of \cite[Theorem 3.2]{tv-cond} proceeds as follows.
\begin{itemize}
\item In Sections 7 and 8, any appearance of $n^{-(2A+1)\gamma}$ should be replaced with $n^{-(2A+2)\gamma+1/2}$.
\item In Section 9, $\eps$ should now be $\frac{1}{\beta} n^{-(2A+1)\gamma+1/2}$ instead of $\frac{1}{\beta} n^{-2A\gamma}$.  Similarly, in (21), the right-hand side should now have an additional factor of $n^{-\gamma+1/2}$.  In Lemma 9.3, the event $X^*_n \in S_{\eps_J,\rho_J}$ should now be $\{ X^*_n \in S_{\eps_J,\rho_J} \} \wedge G$.
\item In Section 11, the problematic sentence ``If we let $u \in {\mathcal N}$ be within $O(\frac{\eps_J}{\rho_J})$ of $X_n^*$, then we have $|X_i \cdot u| \ll \frac{\eps_J}{\rho_J}$...'' should be replaced by ``If we let $u \in {\mathcal N}$ be within $O(\frac{\eps_J}{\rho_J})$ of $X_n^*$, and $G$ holds, then
$$ (\sum_{i=1}^{n-1} |X_i \cdot u|^2)^{1/2} \ll n^\gamma \frac{\eps_J}{\rho_J} $$
by the arguments of the previous section.''  In the next display, the right-hand side should now be
$$ \sum_{u \in {\mathcal N}} {\bf P}( (\sum_{i=1}^{n-1} |X_i \cdot u|^2)^{1/2} \ll n^\gamma \frac{\eps_J}{\rho_J} ).$$
In the next display, an additional factor of $n^{-\gamma+1/2}$ should now be present in the middle and right hand sides.  The phrase ``by Corollary 6.3 (or by a minor modification of Theorem 6.4)'' should now read ``by Theorem 6.4''.
\end{itemize}

The weakening of \cite[Theorem 3.2]{tv-cond} causes similar weakenings of the corollaries mentioned in that paper:
\begin{itemize}
\item The expressions $n^{-(2A+1)\gamma}$ and $n^{(2A+2)\gamma}$ in \cite[Corollary 3.3]{tv-cond} should now be $n^{-(2A+2)\gamma+1/2}$ and $n^{(2A+3)\gamma-1/2}$ respectively.
\item In \cite[Corollary 3.5]{tv-cond}, $(n^{1/2}+\|M\|)^{-2A-1}$ should now be $(n^{1/2}+\|M\|)^{-2A-2} n^{1/2}$.  But \cite[Corollary 3.6]{tv-cond} remains unchanged.
\item In \cite[Corollary 3.7]{tv-cond}, $(n^{1/2}+\|M\|)^{-1+o(1)}$ should now be $(n^{1/2}+\|M\|)^{-2+o(1)} n^{1/2}$.
\item In \cite[Corollary 3.8]{tv-cond}, $n^{-1-\frac{5}{2}A-A^2}$ should now be $n^{-\frac{3}{2}-3A-A^2}$, and $n^{-1-o(1)}$ should now be $n^{-\frac{3}{2}-o(1)}$.
\end{itemize}

We also take this opportunity to record some other minor typos in \cite{tv-cond}:

\begin{itemize}
\item In Theorem 3.4, ``fourth moment bounded $B$'' should be ``fourth moment bounded by $B$''.
\item In Section 4, the display ${\bf P}(s_n(M_n) \ll n^{1/4} L^{-1/2}) \gg n^{-1/2}$ should be deleted and replaced by a period.
\end{itemize}